\pdfoutput=1
\RequirePackage{ifpdf}
\ifpdf 
\documentclass[pdftex]{sigma}
\else
\documentclass{sigma}
\fi

\newtheorem{Theorem}{Theorem}[section]
\newtheorem{Proposition}[Theorem]{Proposition}
\newtheorem{Lemma}[Theorem]{Lemma}

{ \theoremstyle{definition}
\newtheorem{Remark}[Theorem]{Remark}}

\DeclareMathOperator{\vol}{vol}

\begin{document}

\allowdisplaybreaks

\renewcommand{\thefootnote}{$\star$}

\renewcommand{\PaperNumber}{042}

\FirstPageHeading

\ShortArticleName{Twistors and Bi-Hermitian Surfaces of Non-K\"ahler Type}

\ArticleName{Twistors and Bi-Hermitian Surfaces\\
of Non-K\"ahler Type\footnote{This paper is a~contribution to the Special Issue on Progress in Twistor Theory.
The full collection is available at
\href{http://www.emis.de/journals/SIGMA/twistors.html}{http://www.emis.de/journals/SIGMA/twistors.html}}}

\Author{Akira FUJIKI~$^\dag$ and Massimiliano PONTECORVO~$^\ddag$}

\AuthorNameForHeading{A.~Fujiki and M.~Pontecorvo}

\Address{$^\dag$~Research Institute for Mathematical Sciences, Kyoto University, Japan}
\EmailD{\href{mailto:fujiki@math.sci.osaka-u.ac.jp}{fujiki@math.sci.osaka-u.ac.jp}}

\Address{$^\ddag$~Dipartimento di Matematica e Fisica, Universit\`a Roma Tre., Italy}
\EmailD{\href{mailto:max@mat.uniroma3.it}{max@mat.uniroma3.it}}

\ArticleDates{Received November 20, 2013, in f\/inal form April 04, 2014; Published online April 08, 2014}

\Abstract{The aim of this work is to give a~twistor presentation of recent results about bi-Hermitian metrics on compact
complex surfaces with odd f\/irst Betti number.}


\Classification{53C15; 53C28}

\renewcommand{\thefootnote}{\arabic{footnote}}
\setcounter{footnote}{0}

\section{Introduction}

We treat in this work bi-Hermitian surfaces, by which we mean orientable Riemannian conformal four-manifolds $(M,[g])$
with a~pair of $[g]$-orthogonal complex structures $J_+$ and $J_-$ which induce the \textit{same} orientation and are
linearly independent, i.e., there is some $p\in M$ where $J_+(p) \neq\pm J_-(p)$.
The original motivation is given by a~general interest for Riemannian manifolds admitting more than one orthogonal
complex structure~\cite{gss,po97, sa94,sv}.

The case when $J_\pm$ are opposite oriented is also interesting and gives rise to a~holomorphic splitting of the tangent
bundle; see~\cite{ag} for a~complete account.

In what follows it will be important to look at possible relations between the two Lee forms (denoted by $\beta_\pm$) of
the Hermitian metrics $(g,J_\pm)$.
For a~general Hermitian surface $(M,g,J)$ with fundamental $(1,1)$-form $\omega(\cdot,\cdot):=g(\cdot,J\cdot)$ the Lee
form is the $1$-form uniquely def\/ined by the relation $d\omega = \beta\wedge\omega$.
Recall that we are in real dimension four.
The 1-form $\beta$ captures most of the conformal features of the Hermitian metric: $(g,J)$ is K\"ahler if and only if
$\beta=0$, conformally K\"ahler if and only if $\beta$ is exact and is locally conformally K\"ahler, abbreviated by
l.c.K., if and only if $\beta$ is closed.

We will always assume that $M$ is compact and connected; let us recall at this point an important result of
Gauduchon~\cite{ga84} which says in particular that on any compact conformal Hermitian surface $(M,[g],J)$ there always
is a~unique metric (up to homothety) in the conformal class $[g]$ for which the Lee form is co-closed; this metric is
usually called the {\it Gauduchon metric} of $(M,[g],J)$.
In particular one can take $\beta$ harmonic if $(M,[g],J)$ is compact l.c.K.

Our main interest is to study the complex structure of the two surfaces $S_+:=(M,J_+)$ and $S_-:=(M,J_-)$ under the
condition that the f\/irst Betti number $b_1(M)$ is \textit{odd}. In other words we will be concerned with non-K\"ahler
surfaces and it is therefore important to f\/ind complex curves on them.
To this end let
\begin{gather*}
T:=\{p\in M\, |\, J_+(p)=\pm J_-(p)\}
\end{gather*}
denote the set of points where the two complex structures agree up to sign; it plays an important role in understanding
the complex geometry of the surfaces because it turns out to always be the zero-set of a~holomorphic section of a~line bundle:

\begin{Proposition}[\cite{agg, po97}]
The $($possibly empty$)$ set $T$ is the support of a~common effective divisor $\mathbf{T}
\geq 0$ in both surfaces~$S_\pm$.
\end{Proposition}

Before discussing how $T$ can be described as a~complex curve inside the twistor space of $(M,[g])$ and compute the
f\/irst Chern class $c_1(\mathbf{T})$ using the twistor picture, we have the following

\begin{Remark}
From the dif\/ferential geometric viewpoint, $T$ is the closure of the union of all smooth surfaces in $M$ which are
simultaneously $J_\pm$-holomorphic.

To see this, let $C\subset (M,[g],J_\pm)$ be a~smooth surface (real $2$-dimensional) which is both a~$J_+$-holomorphic
and $J_-$-holomorphic curve and consider the tangent space $T_pC$ at any point $p\in C$.
Let $V$ denote this real $2$-dimensional subspace of $T_pM$.
Because $J_\pm$ are assumed to be $g$-compatible the orthogonal complement $W$ is also a~$J_\pm$-complex subspace of
$T_p M$; in each $2$-dimensional subspace $V$ and $W$ the two complex structures $J_\pm$ are rotation by $90^{\circ}$
and induce the same orientation on~$T_p M$ if and only if they are linearly dependent a~$p$.

The complex curve $T$ turns out in general to be non-smooth with several irreducible components meeting transversally,
the smooth open set of each component is a~$J_{\pm}$ holomorphic curve~$C$ as above and~$T$ is the closure of the union
of all such~$C$.

This observation is relevant to the blow-up construction of Cavalcanti--Gualtieri~\cite{cg} because in their case the
exceptional divisor $E$ does not belong to the anti-canonical divisor $-K$, which is supported precisely on $T$.
This shows that $E$ is a~$J_+$-holomorphic curve which is not $J_-$-holomorphic.
\end{Remark}

In order to geometrically describe how $T$ turns out to be a~complex curve with a~natural structure of a~divisor in each
surface $(M,J_\pm)$ we now present a~\textit{twistor approach} to bi-Hermitian metrics in four-dimension~--
following~\cite{po97}~-- which produced new examples with $b_1$ odd~\cite{fp} as we shall discuss later.

An orthogonal almost complex structure is the same as a~smooth section $J:M\to Z$ of the twistor space $Z$ which is the
f\/iber bundle of all linear complex structures at $T_p M$ compatible with the metric and orientation; in
four-dimension, the f\/iber at $p\in M$ is the homogeneous space ${\rm SO}(4)/{\rm U}(2)\cong \mathbb{CP}_1$.
It is known that $Z$ is an almost complex $6$-manifold~-- let $\mathbb{J}$ denote its almost complex structure~-- which
only depends on the f\/ixed conformal structure $[g]$ and orientation of $M$.
The integrability of $J$ is equivalent to the fact that $J:(M,J)\to (Z,\Bbb J)$ is an almost holomorphic map and that
its image $J(M)=:S$ is an almost complex submanifold of $(Z,\mathbb{J})$ which is therefore (tautologically)
biholomorphic to the original complex surface: $S\cong (M,J)$, see for example~\cite{amp}.
Notice that we make no assumption on the integrability of the twistor almost complex structure $\mathbb{J}$.

The twistor space $Z$ is also equipped with a~natural \textit{real structure}, namely the anti-holo\-mor\-phic involution
$\sigma:Z\to Z$ that sends $J\mapsto-J$; objects that are $\sigma$-invariant are then called `real'.
For example, setting $X:=S \amalg \sigma(S)$ we get a~`real' submanifold of $Z$ and as $M$ is compact we can consider
its Poincar\'e dual $X^* \in H^2_{\sigma}(Z,\mathbb{Z})$ which is a~`real' $2$-cohomology class, i.e.\
$\sigma$-invariant.
By Leray--Hirsch theorem this space is $1$-dimensional generated, over $\mathbb{Q}$, by the f\/irst Chern class~$c_1(Z)$
so that $X^* =q\cdot c_1(Z)$ for some $q\in\mathbb{Q}$.

Each f\/iber $L\cong {\rm SO}(4)/{\rm U}(2)\cong \mathbb{CP}_1$ is also $\sigma$-invariant, usually called a~`real twistor line',
and turns out to have normal bundle $\nu_{L/Z}\cong {\mathcal O}_{\mathbb{CP}_{1}}(1) \oplus {\mathcal
O}_{\mathbb{CP}_{1}}(1)$.
Therefore, $c_1(Z)|_{{L}} =4$ and because for each $L$ the two almost complex submanifolds $X$ and $L$ intersect in
exactly two points we conclude that $q=\frac 12$.

This can be used to compute the f\/irst Cherm class of a~Hermitian conformal surface $(M,[g],J)$ from its twistor space
$Z$ by adjunction formula to the almost complex submanifold $S\subset Z$.
Using that $S\cong (M,J)$ and $X=S \amalg \sigma(S)$ we get from adjunction that
\begin{gather}
\label{adj}
c_1(M,J)=[c_1(Z)- X^*]|_{{S}}=\frac 12 c_1(Z)|_{{S}}.
\end{gather}

In our case we have a~bi-Hermitian structure $([g],J_\pm)$ from which we get two complex surfaces in the twistor space
$Z$.
We will use the following notation $S_\pm:=J_\pm(M)$ and similarly $X_\pm:=S_\pm\amalg \sigma(S_\pm)$.
Notice again that there is a~tautological biholomorphism $(M,J_\pm)\cong S_\pm$.

In this notation the set of points $T$ can be identif\/ied with either of the following two subsets in the twistor
space: $X_+\cap S_-$ or $S_+\cap X_-\subset Z$; this exhibits $T$ as an almost complex subvariety of either $S_+$ and
$S_-$; $T$ is therefore a~complex curve in each of the two smooth surfaces, in particular closed in the analytic Zariski
topology and nowhere dense~\cite[Proposition~1.3]{po97}.

\begin{Proposition}
The complex curve $T$ has a~natural structure of divisor $\mathbf{T}$ in both surfa\-ces~$S_\pm$ given by the property
\begin{gather*}
c_1(\mathbf{T})= c_1(S_\pm).
\end{gather*}
In other words, for each compact complex surface $(S_\pm)$ there is a~holomorphic line bundle $F_\pm$ such that
\begin{enumerate}\itemsep=0pt
\item[$(i)$] $c_1(F_\pm)=0$, and
\item[$(ii)$] $\mathbf{T}=F_\pm-K_\pm$ where $K_\pm$ is the canonical line bundle of the surface $S_\pm$.
\end{enumerate}
\end{Proposition}
\begin{proof}
In the twistor space $Z$ we have that $\mathbf{T}\cong X_-\cap S_+\cong X_+\cap S_-$ therefore $c_1(\mathbf{T})=\frac 12
c_1(Z)|_{{S_+}}= \frac 12 c_1(Z)|_{{S_-}}$ and the conclusion follows from~\eqref{adj}.
\end{proof}

\begin{Remark}
In general, the complex curve $T$ has several irreducible components and it gets the structure of ef\/fective
divisor $\mathbf{T}$, with multiplicities, by taking the (only) linear combination such that
$c_1(\mathbf{T})=c_1(J_-)=c_1(J_+)$.
As usual, in what follows $\mathbf{T}$ will denote the divisor as well as the holomorphic line bundle.
\end{Remark}

Although the two complex structures $J_+$ and $J_-$ are in general not biholomorphic, it turns out that they share
common properties.
For this reason we use the notation $J_\pm$ to denote either of them.

As shown in~\cite{ap01}, the f\/lat line bundles $F_+$ and $F_-$ are also isomorphic because they come from the same
representation of $\pi_1(M)$.
In fact, it is important to recall that as f\/lat line bundles on $M$ they correspond to the 1-form $-\frac 12(\beta_+
+\beta_-)$~\cite{ap01} the opposite average of the two Lee forms of $(g,J_\pm)$); by~\cite[p.~420]{agg} this turns out to
be always a~closed 1-form on a~compact $M$.
We will then denote by $F$ the line bundle of zero Chern class $F_\pm$ and we will be looking at the fundamental
equation
\begin{gather}
\label{fund}
\mathbf{T} = F-K,
\end{gather}
which on each of the two surfaces $(M,J_\pm)$ relates the divisor $\mathbf{T}$ to the canonical line bundle $K$.

We will use the fundamental equation~\eqref{fund} above to understand the complex structures of a~bi-Hermitian surface
as well as its Riemannian properties.
Denote, as usual, by $\beta_{\pm}$ the Lee forms of the bi-Hermitian surface $(M,[g],J_\pm)$; we next want to consider
the following conformally invariant conditions:
\begin{enumerate}\itemsep=0pt
\item[1)]
$\beta_+ -\beta_-=0$;
\item[2)]
$\beta_+ +\beta_-=df$;
\item[3)]
$\mathbf{T}=0$.
\end{enumerate}
The f\/irst equality is the \textit{hyper-Hermitian} condition meaning that $J_+$ and $J_-$ span an $S^2$-worth of
complex structures on $M$ like in the hyper-K\"ahler case.
In this situation it was shown by Boyer~\cite{bo88} that $(M,J_\pm)$ must be a~Hopf surface, when $b_1(M)$ is assumed to
be odd.
In this article we will not be concerned with the hyper-Hermitian case and consider a~bi-Hermitian surface to have
exactly two complex structures (up to sign).

The second condition is equivalent to say that $J_+$ and $J_-$ have the same Gauduchon metric in the conformal class
$[g]$ and that the sum of Lee forms vanishes $\beta_+ +\beta_-=0$, for this metric.

Equivalently, the f\/lat line bundle $F$ is holomorphically trivial or in other words: $\mathbf{T}=-K$.
Indeed, by Gualtieri~\cite{gu10} generalized K\"ahler structures are in bijective correspondence with bi-Hermitian
structures satisfying $\mathbf{T}=-K$.

Finally, the third condition says that $J_+$ and $J_-$ are dif\/ferent (up to sign) at {\it every} point; these are
called {\it strongly bi-Hermitian} metrics.
The equation $\mathbf{T}=0$ implies $c_1(S_\pm)=0$ in $H^2(M,\Bbb Z)$ which, as $b_1(M)$ is odd and ${\rm Kod}(S_\pm)\leq 0$,
implies that $S_\pm$ is a~Kodaira, a~Hopf or a~Bombieri--Inoue surface.
However it is shown in~\cite{ad} that only the second possibility can actually occur.
We will explain this point in the next section.

From now on we will indicate by $S$ any of the two surfaces $S_\pm$.
The fundamental equation~\eqref{fund} says that $\mathbf{T}$ is a~numerically anti-canonical divisor in $S$ and it turns
out that $T$ can have at most two connected components.
Setting $t:=b_0(T)$ we obviously have that $t=0$ corresponds to the strongly bi-Hermitian situation, while $t=2$ is
equivalent to the equation $\mathbf{T}=-K$; this has been observed in~\cite{abd} and will be explained in the last section.
In Section~\ref{Section3} we construct a~twistor example with $t=1$ and illustrate a~very interesting,
more general result by~\cite{abd}.

We conclude this introduction with a~short outline of the rest of the paper.

In Section~\ref{Section2} we present some preliminary results from~\cite{ap01,ad,agg, po97} and a~technical lemma.
We then explain a~result of Apostolov and Dloussky, which asserts that the minimal model of a~bi-Hermitian $S$ with
$b_1(S)$ odd can only be a~Hopf surface or a~Kato surface.

In Section~\ref{Section3} we discuss the case of Hopf surfaces which admit bi-Hermitian metrics of all kinds, i.e.\
with all possible values of $t=0,1,2$.
The f\/irst two cases are due to Apostolov--Dloussky and Apostolov--Bailey--Dloussky.

Section~\ref{Section4} is devoted to study Kato surfaces.
We divide them into two main types: Kato surfaces with branches and without.
The second type is best known and we have a~description of all possible anti-canonical divisors on them.
This gives important information concerning existence of bi-Hermitian metrics.
For example: there is no bi-Hermitian metric whatsoever if the minimal model is a~general Enoki or half Inoue surface.
On the positive side, the f\/irst bi-Hermitian metrics on Kato surfaces were constructed explicitly by
LeBrun~\cite{le91} on parabolic Inoue surfaces; we present in this section a~brief outline of a~twistor construction due
to Fujiki--Pontecorvo of bi-Hermitian metrics on hyperbolic and parabolic Inoue surfaces~\cite{fp}.
Because these metrics are anti-self-dual they satisfy $t=2$~\cite{po97}.

Finally, we show that the situation is fairly satisfactory for blown-up hyperbolic Inoue surfaces because they do admit
bi-Hermitian metrics whenever they can; furthermore $t=2$ always in this case.
The situation for {\it intermediate} Kato surfaces, i.e.\
with branches, is still open: there are no examples of bi-Hermitian structures, as far as we know.

\section{Preliminary results}\label{Section2}

The main tool for our study of bi-Hermitian surfaces of non-K\"ahler type is the fundamental equation~\eqref{fund}; we
start this section with the following useful result of Apostolov about the fundamental line bundle $F$:

\begin{Proposition}[\protect{\cite[Lemma~1]{ap01}, \cite[Proposition~2]{ad}}]
Let $(M,[g],J_\pm)$ be any compact bi-Hermi\-tian surface.
Then the flat line bundle $F$ comes from a~real representation of the funda\-men\-tal group~$\pi_1(M)$.
Furthermore, with respect to any $J_+$-standard metric, we have
\begin{gather*}
\deg(F)=-\frac{1}{8\pi}\int_M\|\beta_+ +\beta_-\|^2.
\end{gather*}
In particular, the degree of~$F$ is non-positive and $\deg(F)=0$ if and only if~$F$ is holomorphically trivial if and
only if the sum of Lee forms $\beta_+ +\beta_-$ is an exact $1$-form for $($one and hence$)$ any metric in the conformal
class.
\end{Proposition}

\begin{Remark}
Notice that in general, the degree of a~holomorphic line bundle $L$ computes the volume (with sign) of a~virtual
meromorphic section with respect to a~f\/ixed Gauduchon metric, i.e.\
with $\partial\bar{\partial}$-closed fundamental $(1,1)$-form, otherwise the degree would not be well def\/ined.

In the special case $c_1(L)=0$, the \textit{sign} of $\deg(L)$ is independent of the chosen Gauduchon metric by a~general
result~\cite{lt} and therefore only depends on the representation of the fundamental group~$\pi_1(M)$.
\end{Remark}

In order to deal with non-minimal surfaces we now introduce the following notation: we let $S_0$ denote the minimal
model of any of $S:=(M,J_\pm)$ with blow-down map $\mathbf{b}:S\to S_0$.
Then, $F_0$ will denote the unique real f\/lat line bundle on $S_0$ such that $F=\mathbf{b}^*(F_0)$, because $\pi_1(S) =
\pi_1(S_0)$; while $\mathbf{b}_*:\mathrm{Div}(S) \to \mathrm{Div}(S_0)$ will be the natural projection and $K_0$ the
canonical bundle on $S_0$.
The following simple observation, see also~\cite{dl06}, will be repeatedly used:

\begin{Lemma}
\label{bd}
Let $(M,[g],J_\pm)$ be a~bi-Hermitian surface and let $S_0$ denote the minimal model of any of $S:=(M,J_\pm)$.
Using the notations introduced above, let $\mathbf{T}_0:=\mathbf{b}_*(\mathbf{T})$.
Then the following hold:
\begin{enumerate}\itemsep=0pt
\item[$i)$] $\mathbf{T}_0 = F_0-K_0$ is an effective divisor on $S_0$;
\item[$ii)$] $S$ is obtained by blowing up $S_0$ at points lying on $\mathbf{T}_0$, in particular $t=b_0(\mathbf{T})=b_0(\mathbf{T}_0)$;
\item[$iii)$] $F_0=0$ if and only if $F=0$, otherwise $F_0$ has negative degree on $S_0$.
\end{enumerate}
\end{Lemma}

\begin{proof}
This is a~standard argument~\cite{ap01,dl06,fp, po97}.
For simplicity we prove it for a~one point blow-up $S$ of $S_0$ with exceptional divisor $E$.
In this case the adjunction formula reads $-K=-\mathbf{b}^*K_0-E$~\cite[p.~187]{gh}.
We can add $F$ to both sides of the equation to get $F-K = F-\mathbf{b}^*K_0-E$; taking $\mathbf{b}_*$ gives the
f\/irst statement~i).

Next we compare the total transform $\mathbf{b}^*\mathbf{T}_0$ with the proper transform $\tilde{\mathbf{T}_0}$.
By adjunction formula again, $\mathbf{b}^*\mathbf{T}_0 = \mathbf{b}^*F_0-b^*K_0 = F-K+E$; on the other hand, we always
have $\mathbf{b}^*\mathbf{T}_0 = \tilde{\mathbf{T}_0} +mE$ where $m$ is the multiplicity of the blown-up point $p$ along
$T_0$; we conclude that $\mathbf{T}= \tilde{\mathbf{T}_0} +(m-1)E$.
This shows that $\mathbf{T}$ is ef\/fective only when $m\geq 1$, i.e.\
$p\in \mathbf{T}_0$, and therefore they have the same number of connected components.

It remains to prove part iii) of the statement.
First of all, it is clear that $F_0$ is torsion if and only if $F$ is torsion which however implies $F=0$ and therefore
$F_0=0$.

On the other hand, since $S$ admits a~bi-Hermitian metric, it will be shown below and independently of this lemma that
$S$ and therefore $S_0$ are class-VII surfaces.
In this case the degree map $H^1(\mathbb{R}^+)(\cong \mathbb{R}^+)\to\mathbb{R}$ is isomorphic, modulo torsion,~\cite[Proposition~3.1.13]{lt}
so that for both $F$ and $F_0$ they are trivial if and only if their degree is zero.
We are left to show $F\neq 0$ implies $\deg(F_0)<0$.

For this purpose, let $\omega$ be the K\"ahler form of a~Gauduchon metric on the minimal model $S_0$.
As is well-known~\cite[p.~186]{gh} there exist small neighborhoods $U\subset V$ of $E$ in $S$ such that the line bundle
$[E]$ admits a~Hermitian metric whose Chern form $\rho$ has support in~$V$, is semi-positive in~$U$ and is positive
def\/inite when restricted to each tangent space of~$E$.
This immediately implies that $\tilde{\omega}_\epsilon:= b^* \omega + \epsilon\rho$ is everywhere positive def\/inite
for any small constant $\epsilon$ and therefore is the K\"ahler form of a~Gauduchon metric on~$S$ with respect to which
we can compute degrees
\begin{gather*}
\deg(F) = \int_S c_1(F)\wedge \tilde{\omega}_\epsilon = \int_S c_1(F)\wedge \mathbf{b}^* \omega + \epsilon
\int_{S} c_1(F)\wedge \rho = \deg(F_0) + \epsilon \int_{S} c_1(F)\wedge \rho.
\end{gather*}
Therefore, if by contradiction $\deg(F_0) >0$ we can f\/ind $\epsilon$ small such that $\deg(F)>0$ which
is absurd.
\end{proof}

Notice that, the divisor $\mathbf{T}_0$ on the minimal model may not, in general, come from a~bi-Hermitian metric on
$S_0$.
It will also be shown that the number of its connected components can only be $0$, $1$ or $2$ and that all these
possibilities actually occur.

We now present the following important result; the f\/irst alternative was proved in~\cite{ap01,ad} and the
second one in~\cite{dl06}.

\begin{Proposition}
Let $(M,[g],J_\pm)$ be a~bi-Hermitian surface with $b_1(M)=odd$ and let $S_0$ denote the minimal model of any of
$S:=(M,J_\pm)$.
Then there are two possibilities:
\begin{enumerate}\itemsep=0pt
\item[$i)$] $b_2(S_0)=0$ and $S_0$ is a~Hopf surface, or else
\item[$ii)$] $b_2(S_0)>0$ and $S_0$ is a~Kato surface.
\end{enumerate}
In particular, $M$ is diffeomorphic to $(S^1\times S^3)\# m\overline{\mathbb{CP}}_2$ with $m=b_2(M)$; or is
a~finite quotient of $S^1\times S^3$.
\end{Proposition}

\begin{proof}
We start by showing that the Kodaira dimension must be negative~\cite{ap01}.
Suppose not, from the fundamental equation~\eqref{fund} we have $\deg(K)=\deg(F)-\deg(\mathbf{T})\leq 0$ because
$\deg(F)\leq 0$ and $\mathbf{T}\geq 0$.
Therefore it is enough to prove that the degree of $K$ cannot vanish.
In fact $\deg(K)=0$ is equivalent to say that $\deg(F)=0=\mathbf{T}$ and therefore equivalent to $K=0$.
By Kodaira classif\/ication it is now enough to exclude that a~Kodaira surface can admit bi-Hermitian metrics.

For this purpose we present a~slightly dif\/ferent argument from Apostolov's original proof.
Let $\omega_+$ denote the fundamental $(1,1)$-form of any Hermitian structure $(M,g,J_+)$ on a~Kodaira surface $(M,J_+)$
and let $\Omega=\gamma_+ -i\delta_+$ be a~holomorphic section of the trivial canonical bundle with the property that
$\gamma_+\wedge \gamma_+=\delta_+\wedge\delta_+=\omega_+\wedge\omega_+=\vol(g)$.

Then $\overline{\partial}\Omega=0$ implies that the real and imaginary parts $\gamma_+$ and $\delta_+$ are self-dual
closed symplectic forms since $\partial\gamma_+=\partial\delta_+=0\in\Lambda^{3,0}=0$ and $d=\partial +
\overline{\partial}$ by the integrability of $J_+$.

\looseness=-1
Therefore, whenever $(M,J_+)$ is a~Kodaira surface, both $(M,\delta_+)$ and $(M,\gamma_+)$ provide Thur\-ston examples of
a~compact symplectic manifold which is not K\"ahlerian~\cite{th76}.
Recall that the topological invariants are: $b_1(M)=3$ while the Euler characteristic and the signature both vanish.

The Kodaira complex structure is not tamed by the Thur\-ston symplectic structures and it was noticed by
Salamon~\cite{sa94} that the relation among them is that they def\/ine an almost hyper-Hermitian structure
$\{J_+,I_+,K_+\}$ containing one complex and two symplectic structures: in fact~$I_+$ and~$K_+$ are tamed by~$\gamma_+$
and~$\delta_+$, respectively.

Assume now by contradiction that there is another $g-$orthogonal complex structure $J_-$.
With the same procedure as above we produce two more self-dual symplectic forms $\gamma_-$ and $\delta_-$ but because
$b_+ (M)=2$ we would get that they are linear combination with \textit{constant} coef\/f\/icients of~$\gamma_+$ and~$\delta_+$; this however would imply that the \textit{angle function} between~$J_+$ and~$J_-$ is also constant, forcing
$J_+$ and $J_-$ to span a~hyper-Hermitian structure on $M$~\cite[Proposition~2.5]{amp}; this is however impossible because such
structures can only live on Hopf surfaces, by a~result of Boyer~\cite{bo88}.

We have shown so far that $\deg(K)<0$ and this certainly implies ${\rm Kod}(S)=-\infty$~-- i.e.\
$S$ belongs to class VII in Kodaira classif\/ication of complex surfaces.
However, $\deg(K)<0$ also implies by Lemma~\ref{bd} that $\deg(K_0)<0$ on the minimal model $S_0$.
We conclude that $S_0$ cannot be a~Bombieri--Inoue~\cite{bm73,in74} because Teleman~\cite{te062} shows that these
surfaces have canonical bundle of positive degree.
Therefore $b_2(S_0)=0$ implies $S_0$ is a~Hopf surface by a~theorem of Bogomolov~\cite{bg76}, later clarif\/ied
by~\cite{lyz,te94}.

It remains to discuss the case $b_2(S)>0$; by the main result of~\cite{dl06} is enough to show that~$\mathbf{T}_0$ is
a~non-trivial divisor on $S_0$.
In fact, if by contradiction $\mathbf{T}_0=0$ we will have $c_1(K_0)=c_1(F_0)=0$ so that $c_1^2(S_0)=0$ but this
equation implies $b_2(S_0)=0$ on class~VII surfaces because $b_1=1$ and $b_2^+=0$.
\end{proof}

\section{Hopf surfaces}\label{Section3}

A Hopf surface is a~compact complex surface whose universal cover is $\mathbb{C}^2\setminus \{0\}$.
The aim of the section is to show that they admit a~surprising abundance of bi-Hermitian metrics, with all possible
values of $t=0,1,2$.

It was shown in~\cite{po97} that any conformally-f\/lat metric $[g]$ on a~Hopf surface $M$ admits two ortogonal complex
structures $J_+$ and $J_-$.
Therefore $(M,[g],J_{\pm})$ becomes a~bi-Hermitian surface with anti-self-dual metric.
In particular $\mathbf{T}=-K$ consists of two disjoint smooth elliptic curves of multiplicity $1$ so that $t=2$.

For some special conformally f\/lat Hopf surface, $J_+$ belongs to a~hyper-Hermitian structure $\{I_+,J_+,K_+\}$ and in
some of these cases the same holds for $J_-$~-- i.e.\
some very special Hopf surfaces have two hyper-Hermitian structures.
The divisor $\mathbf{T}$ is zero in this case.

More in general, for bi-Hermitian metrics with $t=0$ there is a~complete result of Apostolov--Dloussky which says that
such metrics exist if and only if any of $(M,J_\pm)$ is a~Hopf surface whose canonical bundle comes from a~\textit{real}
representation of the fundamental group~\cite{ad}.
For sometime these were the only examples of bi-Hermitian metrics with $F\neq 0$.

More recently however, the f\/irst examples of bi-Hermitian metrics with $t=1$ have been constructed in~\cite{abd},
again on Hopf surfaces.

The techniques are very interesting and general: the aim is to construct a~bi-Hermitian metric with a~{\it connected}
divisor $\mathbf{T}=F-K$.
The ingredient for doing it is a~l.c.K.\
metric which we can think of as a~twisted K\"ahler metric with values in a~f\/lat line $L$ (the degree of which will
automatically be positive).

Assuming that it is possible to take $L=-F$, one can contract a~holomorphic section of $\mathbf{T}$ with the l.c.K.\
$(1,1)$-form to obtain a~tensor f\/ield which is in fact an inf\/initesimal deformation of the complex structure,
trivial along $\mathbf{T}$.
It is shown in~\cite{abd} that this def\/ines a~true deformation and the deformed complex structure as well as the
original one are both orthogonal with respect to a~Riemannian conformal metric.
This circle of ideas was used in~\cite{go12} to show that a~surface of K\"ahler type is bi-Hermitian if and only if
admits holomorphic anti-canonical sections.

We can now give a~twistor proof of the existence of bi-Hermitian metrics with $t=1$, for some very special Hopf
surfaces.

\begin{Proposition}
Let $S$ be a~Hopf surface which is an elliptic fiber bundle over $\mathbb{CP}_1$.
Then $S$ admits a~bi-Hermitian metric with $t=1$.
\end{Proposition}
\begin{proof}
The usual Vaisman metric on $S$ is l.c.K., conformally f\/lat and hyper-Hermitian because $S\to \mathbb{CP}_1$ is
a~smooth bundle~\cite{po97}.
In order to f\/ind the twisting l.c.K.\
line bundle $L$ we can use its twistor space $Z$ as follows.
The hyper-Hermitian structure naturally def\/ines a~holomorphic projection $Z\to \mathbb{CP}_1$ so that the normal
bundle of the image $X$ of the l.c.K.\
metric in $Z$ is trivial.
On the other hand~-- because the twistor f\/ibration induces an isomorphism of fundamental groups and therefore of
f\/lat line bundles~-- the twisting line bundle $L$ can also be read--of\/f from the twistor space.
In fact, by the main result of~\cite{po92} $X=-\frac12 K_Z-L$, therefore the triviality of the normal bundle
$X|_{{X}}$ implies $2L=-(K_Z)|_{{S}}$ on the Hopf surface.
By adjunction formula we then get $2L=-K_S$.
Now, if $E$ denotes any f\/iber of the elliptic f\/ibration $S\to\mathbb{CP}_1$ we have $-K_S=2E$ so that $L=E$.
The machinery of~\cite{abd} therefore yields a~bi-Hermitian metric with $\mathbf{T}=F-K=-L-K=-E+2E=E$, in particular
$t=1$.
\end{proof}

The above is just the easiest example of the following general result which uses the fact that every Hopf surface is
l.c.K.~\cite{go}.

\begin{Theorem}[\cite{abd}]
Hopf surfaces admits bi-Hermitian metrics with $t=1$.
$\mathbf{T}$ is supported on an elliptic curve which can have arbitrary $($positive$)$ multiplicity.
\end{Theorem}

Furthermore, bi-Hermitian metrics on blown-up Hopf surfaces were constructed by LeBrun, using a~hyperbolic ansatz for
anti-self-dual metrics with circle symmetry~\cite{le91}; by Kim--Ponte\-cor\-vo~\cite{kp}
using twistor methods and more recently by a~general construction of Cavalcanti--Gualtieri~\cite{cg}.

In~\cite{fp} anti-self-dual bi-Hermitian structures with $\mathbf{T}=-K$ are constructed on any blown-up Hopf surfaces,
blown-up hyperbolic Inoue surfaces and blown-up parabolic Inoue surfaces $S_t$ which are obtained via a~small
deformation $(S_t,-K_t)$ of `anti-canonical pairs' $(S,-K)$, where~$S$ is any hyperbolic Inoue surface.

\section{Kato surfaces}\label{Section4}

A minimal non-K\"ahler surface of Kodaira dimension $-\infty$ and positive second Betti number is said to belong to
class VII$_0^+$.
All known examples are so called Kato surfaces which were introduced in~\cite{ka77} and by def\/intion are compact
complex surfaces $S$ admitting a~global spherical shell: there is a~holomorphic embedding $\phi\colon U\rightarrow S$,
where ${U\subset{\Bbb{C}^2\setminus \{0\}}}$ is a~neighborhood of the unit sphere $S^3$, such that ${S\setminus
\phi(U)}$ is connected.

The following statement summarizes some of the main results about Kato surfaces

\begin{Theorem}[\cite{dl06,dot,ka77}]
For a~surface $S$ in class $VII_0^+$ the following conditions are equivalent and they
imply that $S$ is diffeomorphic to $(S^1\times S^3)\#n\overline{\mathbb{CP}}_2$.
\begin{enumerate}\itemsep=0pt
\item[$1.$]
$S$ is a~Kato surface.
\item[$2.$]
$S$ contains $b_2(S)$ rational curves.
\item[$3.$]
$S$ admits a~divisor $D= G-mK$ with $c_1(G)=0$ and $m\in\mathbb{Z} ^+$.
\end{enumerate}
\end{Theorem}

A divisor of the form $D= G-mK$ with $c_1(G)=0$ is called a~NAC (numerically anti-canonical) divisor in the terminology
of Dloussky~\cite{dl06} because its Chern class is a~multiple of the anti-canonical class.
It is known that $D$ is automatically ef\/fective on $S\in \mathrm{VII}_0^+$.
The smallest $m\geq 1$ for which a~NAC divisor exists is called the {\it index} of the Kato surface $S$.

There are recent important results about Kato surfaces going in dif\/ferent directions; concer\-ning their Hermitian
geometry Brunella proved the following strong result
\begin{Theorem}[\cite{br11}]
Every Kato surface admits l.c.K.\ metrics.
\end{Theorem}

But Kato surfaces are important mainly because they are the only known examples of surfaces in class-VII$^0_+$.
A strong conjecture of Nakamura~\cite[Conjecture~5.5]{na90} says that there should be no other examples in this class.
We only point out here that there is recent important progress in this direction by A.~Teleman~\cite{te10,te13}.

Because on any bi-Hermitian surface $\mathbf{T}$ is a~NAC divisor with $m=1$ we are led to study these Kato surfaces.
For doing this we divide them into classes according to the intersection properties of the $b_2$ rational curves.

First of all, each Kato surface has a~cycle $C$ of rational curves; we then consider two broad classes: Kato surfaces
with no branch~-- also called {\it extreme}~-- that is, every rational curve belongs to some cycle $C$; and surfaces
with branches~-- also called {\it intermediate}~-- in this case the maximal curve consists of the union of one cycle
with chains of rational curves intersecting it transversally at a~single point~\cite{na84}.

We will only be concerned about unbranched Kato surfaces which we are now going to describe in more detail.
These surfaces can be divided into two subclasses; for simplicity we consider their characterizations due to
Nakamura~\cite{na84} in terms of the conf\/iguration of complex curves on them, rather than their original def\/initions
due to Inoue~\cite{in77} and Enoki~\cite{e81}.

We start by considering the class of Enoki surfaces, containing parabolic Inoue surfaces.
The characterizing property of these surfaces is that the unique cycle $C$ has self-intersection zero: $C^2=0$, which
(as $b_2\neq 0$) is the same as saying that every irreducible component of $C$ has self-intersection number $-2$.
Because the intersection form is negative def\/inite on any surface of class VII$^+$, it follows that $C=0$ in
$H^2(M,\Bbb Z)$.
We can now state the following important result of Enoki

\begin{Theorem}[\cite{e81}]
A surface $S\in {\rm VII}_0^+$ has a~divisor $D$ of zero self-intersection if and only if $S$ is an Enoki surface
and $D=mC$ is a~multiple of the cycle of rational curves in $S$.
\end{Theorem}

Enoki surfaces are exceptional compactif\/ications of af\/f\/ine line budles in the sense that the complement
$S\setminus C$ is always an af\/f\/ine line bundle over an elliptic curve.
In the special case that this bundle has a~section the Enoki surface $S$ is special because it contains an elliptic
curve $E$; such special Enoki surfaces are called {\it parabolic Inoue} surfaces, it turns out that they are the only
surfaces in class VII$_0^+$ containing one (and in fact only one) elliptic curve $E$; because the anti-canonical bundle
$-K=E+C$ it follows that $E^2=-b_2(S)$.
We will use the following terminology: a~{\it general Enoki} surface is an Enoki surface which is not parabolic Inoue.

The other class of Kato surfaces without branches is that of {\it hyperbolic Inoue} surfaces, these are the only Kato
surfaces with two cycles of rational curves and are also called even Hirzebruch--Inoue surfaces.
In the special case when these cycles are isomorphic~-- in particular their components have the same intersection numbers
in the same cyclic order~-- the hyperbolic Inoue surface has a~f\/ix-point-free involution whose quotient is a~Kato
surface called {\it half Inoue} surface, or also odd Hirzebruch--Inoue surface.
Finally, every unbranched Kato surface falls into one of these classes.

We can now start reporting about bi-Hermitian metrics on them, most of the results below can also be found
in~\cite[Appendix~A]{abd}.

The f\/irst examples are due to LeBrun~\cite{le91} who used his hyperbolic ansatz to produce anti-self-dual Hermitian
metrics with $S^1$-action on parabolic Inoue surfaces.
It was noticed later on that these metrics are actually bi-Hermitian with $t=2$~\cite{po97}.

Let us notice here that the hyperbolic ansatz can only work in this case because by~\cite{po97} the isometric action is
automatically holomorphic but by a~result of Hausen~\cite{ha95} parabolic Inoue surfaces are the only Kato surfaces
which can admit holomorphic $S^1$-action.

The other known examples are again anti-self-dual and came from the twistor construction of Fujiki--Pontecorvo~\cite{fp}
which we are now going to brief\/ly describe.

We started from a~Joyce twistor space, that is the twistor space $Z$ of a~($S^1\times S^1$)-invariant self-dual metric
on the connected sum $m\mathbb{CP}_2$ of complex projective planes~\cite{jo95}.
These twistor spaces were studied by Fujiki~\cite{fu00} who showed that there are $(m+2)$ invariant twistor lines and
each of them is the transverse intersection of a~pair of {\it elemetary} divisors $S_+$, $S_-$~-- this merely means
that their intersection number with a~twistor line equals~1.

As explained in the introduction a~bi-Hermitian metric produces two pairs of {\it disjoint} complex hypersurfaces in
the twistor space $Z$ which turns out to be a~complex $3$-manifold precisely when the metric is
anti-self-dual~\cite{ahs}.
In order to obtain these desired conf\/iguration of hypersurfaces we therefore blow up two invariant twistor lines and
obtain a~smooth complex $3$-fold $\tilde{Z}$ with two pairs of disjoint hypersurfaces.
Even though $\tilde{Z}$ is not a~twistor space anymore, this is the starting point of a~method introduced in~\cite{df}
in order to construct the twistor space of a~connected sum.

Following this method and its relative version~\cite{kp} we now carefully identify the resulting two exceptional
divisors in $\tilde{Z}$ and consider the resulting quotient space $\hat{Z}$ which is a~`singular twistor space'
containing a~Cartier divisor $\hat{S}$ consisting of two pairs of disjoint singular hypersurfaces.
$\hat{S}$~is in fact an anti-canonical divisor in $\hat{Z}$.

The general theory of~\cite{df} and~\cite{kp} then predicts that if there exist smooth and `real' deformations
$(Z_t,S_t)$ of the singular pair $(\hat{Z},\hat{S})$ then $Z_t$ is a~smooth twistor space containing an ef\/fective
anti-canonical divisor $S_t$ having $4$ irreducible components, pairwise disjoint, and each component meets every
twistor line in exactly one point.
In other words any smooth deformation will yield exactly the twistor data of a~bi-Hermitian anti-self-dual structure on
the connected sum.
The existence of a~real structure on the deformation is a~consequence of the general theory.

Because we started by identifying two exceptional divisors in the same $\tilde{Z}$ the resulting $4$-manifold is
actually a~self-connected sum, dif\/feomorphic to $(S^1\times S^3)\# m \overline{\mathbb{CP}}_2$ which is the
dif\/feomorphism type of a~Kato surface, as desired.

The deformation theory of the constructed singular pair $(\hat{Z},\hat{S})$ is governed by the local ${\rm Ext}$
sheaves and global ${\rm Ext}$ groups of $\Omega_{\hat{Z}}(\log \hat{S})$ which is the sheaf of germs of holomorphic $1$-form
having at worst logarithmic poles along~$\hat{S}$.
We then showed existence of smooth deformations of the pair by proving the vanishing of both
${\rm Ext}^2(\Omega_{\hat{Z}}(\log \hat{S}))$ and $H^2(\Theta_{\hat{Z}}(-\log \hat{S}))$, this last sheaf being holomorphic
vector f\/ields on $\hat{Z}$ which are tangent to $\hat{S}$, along $\hat{S}$.

In this construction the various choices of original $(S^1\times S^1)$-action on $m\mathbb{CP}_2$ and of pairs of
invariant twistor lines we started with allow us to obtain every hyperbolic Inoue surface as an irreducibile component
of~$S_t$.
We can f\/inally state our main result

\begin{Theorem}[\cite{fp}]
Every properly blown up hyperbolic Inoue surface admits families of bi-Hermitian anti-self-dual structures.
The same result holds for  some  parabolic Inoue surfaces.
In particular, $t=2$ for all these metrics.
\end{Theorem}

\begin{Remark}
Properly blown up means that we are only allowed to blow up points which are nodes of (the anti-canonical cycle)~$-K$.

Furthermore, a~variation of the construction gives existence of families of anti-self-dual Hermitian metrics on every
properly blown up half Inoue surface.
These surfaces however cannot admit bi-Hermitian metrics as will soon be clear.

The above were the f\/irst examples of l.c.K.\
metrics without symmetries on Kato surfaces.
\end{Remark}

Now that we have an existence result we can describe the situation for bi-Hermitian metrics on unbranched Kato surfaces;
this is made possible by the fact that all their NAC divisors can be described.
A useful consequences of Enoki theorem is that the rational curves on a~Kato surface form a~bases of $H^2(S,\mathbb{Q})$
unless $S$ is Enoki.
In particular a~NAC divisor $D=G-K$ is unique or else the surface is Enoki.

As usual, the statements below hold for any one of the two complex surfaces $S=(M,J_\pm)$.

\begin{Proposition}
There are no bi-Hermitian metrics at all if the minimal model of $S$ is a~general Enoki or half Inoue surface.
\end{Proposition}

\begin{proof}
By Lemma~\ref{bd} it is enough to show that the minimal model contains no NAC divisors of index $1$: $D=G-K$.
We start by considering a~general Enoki surface, it is known in this case that $-K=C+\hat E$ where $C$ is the unique
cycle and $\hat E$ is a~line bundle without meromorphic sections.
The second cohomology $H^2(S,\Bbb Z)$ is spanned by the Chern class $c_1(\hat E)$ and the irreducible components of the
cycle $C$ which however are subject to the relation $C=0\in H^2(S,\Bbb Z)$.
If there was a~NAC divisor $D=G-K$ on a~general Enoki surface we would have $D-C=G+\hat E$.
But this equation is impossible in cohomology because $C$ is the maximal curve so that $c_1(D-C)$ is spanned by the
irreducible components of the cycle while $c_1(G+\hat E)=c_1(\hat E)$ which is linearly independent from the components
of $C$.
This argument shows that there is no NAC divisor on a~general Enoki surface.

Suppose now that the minimal model is a~half Inoue surface.
In this case $-K=C+L$ where~$C$ is the unique cycle and $L$ is a~non-trivial line bundle of order~$2$.
Then the unique NAC divisor is~$C$ which however does not satisfy the fundamental equation~\eqref{fund} $C=F-K$ because
the degree of~$F$ should vanish if and only if~$F$ is the trivial line bundle, while the degree of~$L$ vanishes~-- since
$2L={\mathcal O}_S$~-- even though $L$ is non-trivial.
\end{proof}

We can use this techniques to show the following result stated in the introduction

\begin{Proposition}
Let $M$ be a~compact four manifold with odd first Betti number.
A bi-Her\-mi\-tian surface $S$ satisfies $\mathbf{T}=-K$ if and only if $T$ is disconnected, if and only if $t=2$.
Therefore its minimal model $S_0$ is either $($a finite quotient of$)$ a~diagonal Hopf surface; or a~parabolic Inoue
surface or a~hyperbolic Inoue surface.
\end{Proposition}

\begin{proof}
One direction was proved in~\cite{agg}.
Therefore, suppose that $T$ is disconnected we need to prove $\mathbf{T}=-K$.
By Lemma~\ref{bd} the minimal model $S_0$ has a~disconnected NAC divisor $\mathbf{T}_0=F_0-K_0$ with $\deg(F_0)\leq 0$ and is
enough to show $\deg(F_0)=0$.

First of all, $S_0$ cannot be a~Kato surface with branches because it is shown in~\cite{do99} that any NAC divisor is
supported on the maximal curve which is connected.
Therefore $S_0$ is either a~Hopf or an unbranched Kato surface with a~disconnected NAC divisor $\mathbf{T}_0=F_0-K_0$
and of course $t=2$ otherwise $S_0$ is an elliptic Hopf surface with $\vol(\mathbf{T}_0)$ bigger than $\vol(-K_0)$, which
is impossible.
Suppose~$S_0$ is a~Hopf surface, if it has a~disconnected divisor it must be a~diagonal one in which case $-K_0=E_1+E_2$
is a~union of two elliptic curves with multiplicity one.
Because~$\mathbf{T}_0$ is an ef\/fective and disconnected divisor $\mathbf{T}_0=aE_1 +bE_2$ with $a,b\geq 1$.
The conclusion is that $F_0=(a-1)E_1 +(b-1)E_2$ has positive degree unless $a=b=1$ so that $F_0=0$, as wanted.

We can now assume that $S_0$ is an unbranched Kato surface with ef\/fective and disconnected NAC divisor $T_0=F_0-K_0$.
Now, a~general Enoki has no NAC divisors, the only NAC divisor of a~half Inoue surface is the cycle $C$ which is
connected and so we are left with parabolic and hyperbolic Inoue surfaces which both have ef\/fective $-K_0$; but in
this case $F_0=T_0 + K_0$ is a~divisor and satisf\/ies the hypothesis of Enoki theorem so that $S_0$ is
a~parabolic Inoue surface with $F_0=nC$ leading to $T_0=nC+C+E=(n+1)C+E$.
But $n<0$ because $F_0$ has negative degree therefore $T_0$ ef\/fective implies $n=-1$ in which case $T_0=E$ is
connected.
\end{proof}

We now come back to existence of bi-Hermitian metrics and consider hyperbolic Inoue surfaces; in this case the situation
is fairly clear (and best possible):

\begin{Proposition}
Let $S_0$ be an arbitrary hyperbolic Inoue surface and suppose $S$ is a~blow up of~$S_0$.
Then every bi-Hermitian metric on $S$ satisfies $t=2$.
Furthermore, such metrics exist if and only if $S$ is obtained by blowing up points on the anti-canonical divisor of~$S_0$.
\end{Proposition}

\begin{proof}
The anti-canonical divisor of $S_0$ is ef\/fective: in fact $-K=C_1 + C_2$ is the maximal curve consisting of the union
of two cycles.
This is the only NAC divisor with $m=1$ by uniqueness and using Lemma~\ref{bd} this proves the f\/irst part of the statement.

To prove the second part, recall that anti-self-dual bi-Hermitian metrics were constructed in~\cite{fp} on every
properly blown up hyperbolic Inoue surface.
Therefore using the result of~\cite{cg}~-- which says that generalizad K\"ahler metrics persist by blowing up
{\it smooth}
points of $-K$~-- we get such metrics on any $S$ obtained by blowing up points on the anti-canonical divisor of its
minimal model.

It is important to notice that we get all possible blow ups because the anticanonical divisor of a~hyperbolic Inoue
surface is reduced; therefore the exceptional curve obtained by blowing up a~smooth point will not belong to the
anticanonical divisor of the blown up surface.
\end{proof}

A similar but weaker existence result holds for bi-Hermitian metrics with $T=-K$ on blown up parabolic Inoue surfaces,
with an additional reality condition.
It remains to be seen whether there are bi-Hermitian metrics with $t=1$ on parabolic Inoue surfaces and on intermediate
(i.e.~branched) Kato surfaces.

\subsection*{Acknowledgements}

We thank the referee for several suggestions which led to signif\/icant improvements in the exposition.

\pdfbookmark[1]{References}{ref}
\LastPageEnding

\end{document}